\documentclass[12pt, a4j]{article}
\AtBeginDvi{}
\usepackage{ascmac}
\usepackage{amsmath, amssymb, amscd}
\usepackage{type1cm}
\usepackage{bm}
\usepackage[pdftex,hiresbb]{graphicx}
\usepackage{float}
\usepackage{tikz}
\newtheorem{defi}{Definition}[section]
\newtheorem{rema}[defi]{Remark}
\newtheorem{exa}[defi]{Example}
\newtheorem{propo}[defi]{Proposition}
\newtheorem{thm}[defi]{Theorem}

\newtheorem{cor}[defi]{Corollary}
\newtheorem{lem}[defi]{Lemma}
\makeatletter
    
    \@addtoreset{equation}{section}
  \makeatother

\begin{document}
\begin{center}
\LARGE \textbf{On volume functions of special flow polytopes} \\
\vspace{0.5\baselineskip}
\normalsize Takayuki NEGISHI, Yuki SUGIYAMA, and Tatsuru TAKAKURA
\end{center}

\vspace{0.5cm}

\begin{abstract}
In this paper, we consider the volume of a special kind of flow polytope. 
We show that its volume satisfies a certain system of differential equations, 
and conversely, the solution of the system of differential equations 
is unique up to a constant multiple. 
In addition, we give an inductive formula for the 
volume with respect to the rank of the root system of type A. 
\end{abstract}

\section{Introduction}

The number of lattice points and the volume of a convex polytope are important and 
interesting objects and have been studied from various points of view (see, e.g., \cite{4}). 
For example, the number of lattice points of a convex polytope associated to a root system is 
called the Kostant partition function, and it plays an important role 
in representation theory of Lie groups (see, e.g., \cite{7}). 

In this paper, we consider a convex polytope 
associated to the root system of type $A$, 
which is called a {\it flow polytope}. 
As explained in \cite{2,3}, 
the cone spanned by the positive roots is divided into several 
polyhedral cones called {\it chambers}, and 
the combinatorial property of a flow polytope depends on a chamber. 
Moreover, there is a specific chamber called the {\it nice chamber}, 
which plays a significant role in \cite{9}. 
Also in \cite{2,3}, a number of theoretical results related to 
the Kostant partition function and the volume function of a flow polytope can be found. 
In particular, it is shown that these functions for the nice chamber are 
written as iterated residues (\cite[Lemma 21]{3}). 
We also refer to \cite{1} for similar formulas for other chambers in more general settings. 

The purpose of this paper is to characterize the volume function of a flow polytope 
for the nice chamber in terms of a system of differential equations, 
based on a result in \cite{3}. 
In order to state the main results, 
we give some notation. 
Let $e_i$ be the standard basis of $\mathbb{R}^{r+1}$ and let 
\[ A_{r}^{+} = \{ e_i - e_j \,|\, 1 \leq i < j \leq r+1 \} \] 
be the positive root system 
of type $A$ with rank $r$. We assign a positive integer $m_{i,j}$
to each $i$ and $j$ with $1 \leq i < j \leq r+1$. 
Let us set $m = (m_{i,j})$ and 
$M = \sum_{1 \leq i < j \leq r+1} m_{i,j}$. 
For $a = a_{1}e_1 + \cdots + a_{r}e_r - (a_1 + \cdots + a_r)e_{r+1} 
\in \mathbb{R}^{r+1}$, where $a_i \in \mathbb{R}_{\geq 0} \,\, (i = 1, \dots, r)$, 
the following polytope $P_{A_{r}^{+}, m}(a)$ is called the flow polytope: 
\begin{eqnarray*}
P_{A_{r}^{+}, m}(a) = \left\{ (y_{i,j,k}) \in \mathbb{R}^M \,\middle|\, 
\begin{array}{l} 1 \leq i < j \leq r + 1 \, , \,\, 1 \leq k \leq m_{i,j} \, , \\ 
y_{i,j,k} \geq 0 \, , \,\, \sum_{1 \leq i < j \leq r+1} 
\sum_{1 \leq k \leq m_{i,j}} y_{i,j,k} (e_i - e_j) = a 
\end{array} \right\}. 
\end{eqnarray*}
Note that the flow polytopes in \cite{3} include the case 
that some of $m_{i,j}$'s are zero, whereas we exclude such cases in this paper. 
We denote the volume of $P_{A_{r}^{+}, m}(a)$ by $v_{A_{r}^{+}, m}(a)$.

The open set 
$$
\mathfrak{c}_{\text{nice}}:=\{ a = a_{1}e_1 + \cdots + a_{r}e_r - (a_1 + \cdots + a_r)e_{r+1} \in \mathbb{R}^{r+1} \,|\, a_i > 0, i = 1, \dots, r \}
$$
in ${\mathbb R}^{r+1}$ is called the nice chamber. 
We are interested in the volume $v_{A_{r}^{+}, m}(a)$ 
when $a$ is in the closure 
of the nice chamber, and then it is written by $v_{A_{r}^{+}, m, \mathfrak{c}_{\rm nice}}$. It is a homogeneous polynomial of degree $M - r$.
The first result of this paper is the following. 

\begin{thm} 
Let $a = \sum_{i = 1}^{r} a_{i}(e_i - e_{r+1}) \in \overline{\mathfrak{c}_{{\rm nice}}}$, 
and let $v_{A_{r}^{+}, m, \mathfrak{c}_{\rm nice}}(a)$ be the volume of $P_{A_{r}^{+}, m}(a)$. Then $v = v_{A_{r}^{+}, m, \mathfrak{c}_{\rm nice}}(a)$ satisfies the system of differential equations as follows:
\begin{equation*}
\begin{cases}
\partial_{r}^{m_{r,r+1}} v = 0 \\
(\partial_{r-1} - \partial_{r})^{m_{r-1,r}} \partial_{r-1}^{m_{r-1,r+1}} v = 0 \\
\,\,\,\,\,\,\,\,\,\, \vdots  \\
(\partial_1 - \partial_2)^{m_{1,2}} (\partial_1 - \partial_3)^{m_{1,3}} \cdots (\partial_1 - \partial_r)^{m_{1,r}} \partial_{1}^{m_{1,r+1}} v = 0,
\end{cases}
\end{equation*} 
where $\partial_i = \frac{\partial}{\partial a_i}$ for $i = 1, \dots, r$. Conversely, the polynomial $v = v(a)$ of degree $M - r$ satisfying the above equations is equal to a constant multiple of $v_{A_{r}^{+}, m, \mathfrak{c}_{\rm nice}}(a)$. 
\end{thm}

We remark that it is known the volume function $v_{A_{r}^{+}, m}(a)$ of $P_{A_{r}^{+}, m}(a)$, as a distribution on $\mathbb{R}^{r}$, satisfies the differential equation 
\begin{equation}
Lv_{A_{r}^{+}, m}(a) = \delta(a)
\end{equation}
in general, where $L = \prod_{i < j} (\partial_i - \partial_j)^{m_{i,j}}$ and $\delta(a)$ is the Dirac delta function on $\mathbb{R}^{r}$ (\cite{6,9}). Note that $\partial_{r+1}$ in the definition of $L$ is supposed to be zero. The above theorem characterizes the function $v_{A_{r}^{+}, m, \mathfrak{c}_{\rm nice}}(a)$ on $\overline{\mathfrak{c}_{\rm nice}}$ more explicitly.



In addition, in Theorem 3.6, we show the volume $v_{A_{r}^{+}, m, \mathfrak{c}_{\rm nice}}(a)$ is written by a linear combination of $v_{A_{r-1}^{+}, m', \mathfrak{c}'_{\rm nice}}(a)$ and its partial derivatives, where $m' = (m_{i,j})_{2 \leq i < j \leq r+1}$ and $\mathfrak{c}'_{\rm nice}$ is the nice chamber of $A_{r-1}^{+}$.

This paper is organized as follows. In Section 2, we recall the iterated residue, the Jeffrey-Kirwan residue, and the nice chamber based on \cite{2}, \cite{3}, \cite{5} and \cite{8}. 
Also, 
we give the some examples of $P_{A_{r}^{+}, m}(a)$ and the calculations of the volume $v_{A_{r}^{+}, m, \mathfrak{c}_{\rm nice}}(a)$. 
In Section 3, we prove the main theorems.

\section{Preliminaries}

In this section, we set up the tools to prove the main theorems based on \cite{2}, \cite{3}, \cite{5} and \cite{8}.

\subsection{Flow polytopes and its volumes}

Let $e_1, \dots, e_{r+1}$ be the standard basis of $\mathbb{R}^{r+1}$, and let 
\[ V = \left\{ a = \sum_{i=1}^{r+1} a_{i}e_i \in \mathbb{R}^{r+1} \,\middle|\, \sum_{i=1}^{r+1} a_i = 0 \right\}. \]
We consider the positive root system of type $A$ with rank $r$ as follows:
\[ A_{r}^{+} = \{ e_i - e_j \,|\, 1 \leq i < j \leq r+1 \}. \]

Let $C(A_{r}^{+})$ be the convex cone generated by $A_{r}^{+}$:
\[ C(A_{r}^{+}) = \{ a = a_{1}e_1 + \cdots + a_{r}e_r - (a_1 + \cdots + a_r)e_{r+1} \,|\, a_1, \dots, a_r \in \mathbb{R}_{\geq 0} \}. \]

We assign a positive integer $m_{i,j}$ to each $i$ and $j$ with $1 \leq i < j \leq r + 1$, and it is called a multiplicity. Let us set $m = (m_{i,j})$ and $M = \sum_{1 \leq i < j \leq r+1} m_{i,j}$.


\begin{defi} 
\textup{Let $a = a_{1}e_1 + \cdots + a_{r}e_r - (a_1 + \cdots + a_r)e_{r+1} \in C(A_{r}^{+})$. We consider the following polytope:
\begin{eqnarray*}
P_{A_{r}^{+}, m}(a) = \left\{ (y_{i,j,k}) \in \mathbb{R}^{M} \,\middle|\, \begin{array}{l} 1 \leq i < j \leq r + 1 \, , \,\, 1 \leq k \leq m_{i,j} \, , \\ y_{i,j,k} \geq 0 \, , \,\, \sum_{1 \leq i < j \leq r+1} \sum_{1 \leq k \leq m_{i,j}} y_{i,j,k} (e_i - e_j) = a \end{array} \right\}, 
\end{eqnarray*}
which is called the {\it flow polytope}.}
\end{defi}

\begin{rema} 
\textup{The flow polytopes in \cite{3} include the case that $m_{i,j} = 0$ for some $i$ and $j$.}
\end{rema}

The elements of $A_{r}^{+}$ generate a lattice $V_{\mathbb{Z}}$ in $V$. The lattice $V_{\mathbb{Z}}$ determines a measure $da$ on $V$. 

Let $du$ be the Lebesgue measure on $\mathbb{R}^{M}$. Let $[ \alpha_1, \dots, \alpha_M ]$ be a sequence of elements of $A_{r}^{+}$ with multiplicity $m_{i,j}$, and let $\varphi$ be the surjective linear map from $\mathbb{R}^{M}$ to $V$ defined by $\varphi(e_k) = \alpha_k$. The vector space $\ker (\varphi) = \varphi^{-1}(0)$ is of dimension $d = M - r$ and it is equipped with the quotient Lebesgue measure $du / da$. For $a \in V$, the affine space $\varphi^{-1}(a)$ is parallel to $\ker (\varphi)$, and thus also equipped with the Lebesgue measure $du / da$. Volumes of subsets of $\varphi^{-1}(a)$ are computed for this measure. In particular, we can consider the volume $v_{A_{r}^{+}, m}(a)$ of the polytope $P_{A_{r}^{+}, m}(a)$. 


\subsection{Total residue and iterated residue}

Let $A_{r} = A_{r}^{+} \cup (-A_{r}^{+})$, and let $U$ be the dual vector space of $V$. We denote by $R_{A_{r}}$ the ring of rational functions $f(x_1, \dots, x_r)$ on the complexification $U_{\mathbb{C}}$ of $U$ with poles on the hyperplanes $x_i - x_j = 0 \,\, (1 \leq i < j \leq r+1)$ or $x_i = 0 \,\, (1 \leq i \leq r)$. A subset $\sigma$ of $A_r$ is called a {\it basis} of $A_r$ if the elements $\alpha \in \sigma$ form a basis of $V$. In this case, we set
\[ f_{\sigma}(x) := \frac{1}{\prod_{\alpha \in \sigma} \alpha(x)} \]
and call such a element a {\it simple fraction}. We denote by $S_{A_r}$ the linear subspace of $R_{A_r}$ spanned by simple fractions. The space $U$ acts on $R_{A_r}$ by differentiation: $(\partial(u)f)(x) = (\frac{d}{d\varepsilon})f(x + \varepsilon u) |_{\varepsilon = 0}$. We denote by $\partial(U)R_{A_r}$ the space spanned by derivatives of functions in $R_{A_r}$. It is shown in \cite[Proposition 7]{5} that
\[ R_{A_r} = \partial(U)R_{A_r} \oplus S_{A_r}. \]
The projection map ${\rm Tres}_{A_r} : R_{A_r} \to S_{A_r}$ with respect to this decomposition is called the total residue map.

We extend the definition of the total residue to the space $\hat{R}_{A_r}$ consisting of functions $P / Q$ where $Q$ is a finite product of powers of the linear forms $\alpha \in A_r$ and $P = \sum_{k=0}^{\infty} P_k$ is a formal power series with $P_k$ of degree $k$. As the total residue vanishes outside the homogeneous component of degree $-r$ of $A_r$, we can define ${\rm Tres}_{A_r} (P/Q) = {\rm Tres}_{A_r} (P_{q-r}/Q)$, where $q$ is degree of $Q$. For $a \in V$ and multiplicities $m = (m_{i,j}) \in (\mathbb{Z}_{\geq 0})^{M}$ of elements of $A_{r}^{+}$, the function 
\[ F := \frac{e^{a_{1}x_1 + \cdots + a_{r}x_r}}{\prod_{i=1}^{r} x_{i}^{m_{i,r+1}} \prod_{1 \leq i < j \leq r} (x_i - x_j)^{m_{i,j}}} \]
is in $\hat{R}_{A_r}$. We define $J_{A_{r}^{+}, m}(a) \in S_{A_r}$ by
\[ J_{A_{r}^{+}, m}(a) = {\rm Tres}_{A_r} F. \]

Next, we describe the iterated residue.

\begin{defi} 
\textup{For $f \in R_{A_r}$, we define the iterated residue by
\[ {\rm Ires}_{x = 0} f = {\rm Res}_{x_1 = 0}{\rm Res}_{x_2 = 0} \cdots {\rm Res}_{x_r = 0} f(x_1, \dots, x_r). \]
}
\end{defi}

Since the iterated residue ${\rm Ires}_{x=0} f$ vanishes on the space $\partial(U)R_{A_{r}}$ as in \cite{3}, we have
\begin{equation}
{\rm Ires}_{x=0} J_{A_{r}^{+}, m}(a) = {\rm Ires}_{x=0} F.  
\end{equation}


\subsection{Chambers and Jeffrey--Kirwan residue}

\begin{defi} 
\textup{Let $C(\nu)$ be the closed cone generated by $\nu$ for any subset $\nu$ of $A_{r}^{+}$ and let $C(A_{r}^{+})_{{\rm sing}}$ be the union of the cones $C(\nu)$ where $\nu$ is any subset of $A_{r}^{+}$ of cardinal strictly less than $r = \dim V$. By definition, the set $C(A_{r}^{+})_{{\rm reg}}$ of $A_{r}^{+}$-regular elements is the complement of $C(A_{r}^{+})_{{\rm sing}}$. A connected component of $C(A_{r}^{+})_{{\rm reg}}$ is called a {\it chamber}.}
\end{defi}


The Jeffrey--Kirwan residue \cite{8} associated to a chamber $\mathfrak{c}$ of $C(A_{r}^{+})$ is a linear form $f \mapsto \langle\langle \mathfrak{c}, f \rangle\rangle$ on the vector space $S_{A_{r}}$ of simple fractions. Any function $f$ in $S_{A_{r}}$ can be written as a linear combination of functions $f_{\sigma}$, with a basis $\sigma$ of $A_r$ contained in $A_{r}^{+}$. To determine the linear map $f \mapsto \langle\langle \mathfrak{c}, f \rangle\rangle$, it is enough to determine it on this set of functions $f_{\sigma}$. So we assume that $\sigma$ is a basis of $A_{r}$ contained in $A_{r}^{+}$.


\begin{defi} 
\textup{For a chamber $\mathfrak{c}$ and $f_{\sigma} \in S_{A_r}$, we define the Jeffrey--Kirwan residue $\langle\langle \mathfrak{c}, f_{\sigma} \rangle\rangle$ associated to a chamber $\mathfrak{c}$ as follows:
\begin{itemize}
\item If $\mathfrak{c} \subset C(\sigma)$, then  $\langle\langle \mathfrak{c}, f_{\sigma} \rangle\rangle = 1$.
\item If $\mathfrak{c} \cap C(\sigma) = \emptyset$, then  $\langle\langle \mathfrak{c}, f_{\sigma} \rangle\rangle = 0$,
\end{itemize}
where $C(\sigma)$ is the convex cone generated by $\sigma$.}
\end{defi}

\begin{rema} 
\textup{More generally, as in \cite[Definition 11]{3}, the Jeffrey--Kirwan residue $\langle\langle \mathfrak{c}, f_{\sigma} \rangle\rangle$ is defined to be $\frac{1}{{\rm vol}(\sigma)}$ if $\mathfrak{c} \subset C(\sigma)$, where ${\rm vol}(\sigma)$ is the volume of the parallelepiped $\oplus_{\alpha \in \sigma} [0, 1]\alpha$, relative to our Lebesgue measure $da$. In our case, the volume ${\rm vol}(\sigma)$ is equal to $1$ since $A_r$ is unimodular.} 
\end{rema}

The volume $v_{A_{r}^{+}, m}(a)$ of the flow polytope $P_{A_{r}^{+}, m}(a)$ is written by the function $J_{A_{r}^{+}, m}(a)$ and the Jeffrey--Kirwan residue in the following.

\begin{thm}[\cite{3}] 
Let $\mathfrak{c}$ be a chamber of $C(A_{r}^{+})$. Then, for $a \in \bar{\mathfrak{c}}$, the volume $v_{A_{r}^{+}, m}(a)$ of $P_{A_{r}^{+}, m}(a)$ is given by
\[ v_{A_{r}^{+}, m}(a) = \langle\langle \mathfrak{c}, J_{A_{r}^{+}, m}(a) \rangle\rangle . \]
\end{thm}

We denote by $v_{A_{r}^{+}, m, \mathfrak{c}}(a)$ the polynomial function of $a$ coinciding with $v_{A_{r}^{+}, m}(a)$ when $a \in \bar{\mathfrak{c}}$. It is a homogeneous polynomial of degree $M - r$.

\subsection{Nice chamber}






\begin{defi} 
\textup{The open subset $\mathfrak{c}_{{\rm nice}}$ of $C(A_{r}^{+})$ is defined by
\[ \mathfrak{c}_{{\rm nice}} = \{ a \in C(A_{r}^{+}) \,|\, a_i > 0 \,\, (i = 1, \dots, r) \}. \]
The set $\mathfrak{c}_{{\rm nice}}$ is in fact a chamber for the root system $A_{r}^{+}$ (\cite{3}). The chamber $\mathfrak{c}_{{\rm nice}}$ is called the {\it nice chamber}.}  
\end{defi}

\begin{lem}[\cite{3}] 
For the nice chamber $\mathfrak{c}_{{\rm nice}}$ of $A_{r}^{+}$ and $f \in S_{A_{r}}$, we have
\[ \langle\langle \mathfrak{c}_{{\rm nice}}, f \rangle\rangle = {\rm Ires}_{x = 0} f. \]
\end{lem}

From Theorem 2.7, Lemma 2.9 and (2.1), we have the following corollary.
 
\begin{cor} 
Let $a \in \overline{\mathfrak{c}_{{\rm nice}}}$. Then the volume function $v_{A_{r}^{+}, m, \mathfrak{c}_{{\rm nice}}}(a)$ is given by
\[ v_{A_{r}^{+}, m, \mathfrak{c}_{{\rm nice}}}(a) = {\rm Ires}_{x=0} F. \]
\end{cor}

\subsection{Examples}

In this subsection, we give some examples of the flow polytopes for $A_1, A_2$, and $A_3$, and calculate their volumes.

\begin{exa} 
\textup{When $r = 1$, the nice chamber of $A_{1}^{+}$ is $\mathfrak{c}_{{\rm nice}} = \{ a = a_1(e_1 - e_2) \,|\, a_1 > 0 \}$. For $a = a_1(e_1 - e_2) \in \overline{\mathfrak{c}_{{\rm nice}}}$,}
\[ P_{A_{1}^{+}, m}(a) = \left\{ (y_{i,j,k}) \in \mathbb{R}^{m_{1,2}} \,|\, y_{i,j,k} \geq 0 \, , \,\, y_{1,2,1} + y_{1,2,2} + \cdots + y_{1,2,m_{1,2}} = a_1 \right\}. \]
\textup{From Corollary 2.10, we have}
\begin{align*}
v_{A_{1}^{+}, m, \mathfrak{c}_{\rm nice}}(a) &= {\rm Res}_{x_1 = 0} \left( \frac{e^{a_{1}x_1}}{x_{1}^{m_{1,2}}} \right) \\
&= \frac{1}{(m_{1,2} - 1) !}a_{1}^{m_{1,2} - 1}.
\end{align*}
\end{exa}

\begin{exa} 
\textup{When $r = 2$, there are two chambers $\mathfrak{c}_1, \mathfrak{c}_2$ of $A_{2}^{+}$ as below, and the nice chamber $\mathfrak{c}_{{\rm nice}}$ of $A_{2}^{+}$ is $\mathfrak{c}_{1}$.}

\begin{center}
\begin{tikzpicture}[line width=1.5pt] 
\draw[-stealth] (0,0) -- (5,0) node[below] {$e_1 - e_2$} ;
\draw[-stealth] (0,0) -- (2.5,4.5) node[right] {$e_1 - e_3$} ;
\draw[-stealth] (0,0) -- (-2.5,4.5) node[right] {$e_2 - e_3$} ;
\draw (2.5,1.5) node {$\mathfrak{c}_2$} ;
\draw (0,2.5) node {$\mathfrak{c}_1$} ;
\draw (-0.25,-0.25) node {$\mathbf{O}$} ;
\end{tikzpicture}
\end{center}
\begin{center}
\textup{Figure 1 \,:\, The chamber of $A_{2}^{+}$.}
\end{center}

\textup{For example, we set $m_{1,2} = n \,\, (n \in \mathbb{Z}_{> 0})$, $m_{1,3} = 1$, and $m_{2,3} = 1$. For $a = a_{1}e_1 + a_{2}e_2 - (a_1 + a_2)e_3 \in \overline{\mathfrak{c}_{{\rm nice}}}$,}
\begin{eqnarray*}
P_{A_{2}^{+}, m}(a) = \left\{ (y_{i,j,k}) \in \mathbb{R}^{n+2} \,\middle|\, \begin{array}{ll} y_{i,j,k} \geq 0 \\ y_{1,2,1} + y_{1,2,2} + \cdots + y_{1,2,n} + y_{1,3,1} = a_1 \\ - y_{1,2,1} - y_{1,2,2} - \cdots - y_{1,2,n} + y_{2,3,1} = a_2 \end{array} \right\}. 
\end{eqnarray*}
\textup{From Corollary 2.10, we have}
\begin{align*}
v_{A_{2}^{+}, m, \mathfrak{c}_{\rm nice}}(a) &= {\rm Ires}_{x = 0} \left( \frac{e^{a_{1}x_1 + a_{2}x_2}}{x_{1}x_2(x_1 - x_2)^n} \right) \\
&= {\rm Res}_{x_1 = 0}{\rm Res}_{x_2 = 0} \left( \frac{e^{a_{1}x_1 + a_{2}x_2}}{x_{1}x_2(x_1 - x_2)^n} \right) \\
&= \frac{1}{n !}a_{1}^n.
\end{align*}
\end{exa}

\begin{exa} 
\textup{When $r = 3$, there are seven chambers of $A_{3}^{+}$ as below (\cite{1}), and the nice chamber $\mathfrak{c}_{{\rm nice}}$ of $A_{3}^{+}$ is $\mathfrak{c}_{1}$.}

\begin{center}
\begin{tikzpicture} 
\draw[line width=1.5pt] (0,0) -- (6,0) -- (3,5) -- (0,0) ;
\draw (0,0) -- (4.5,2.5) ;
\draw (6,0) -- (1.5,2.5) node[left] {$e_1 - e_3$} ;
\draw (1.5,2.5) -- (4.5,2.5) node[right] {$e_2 - e_4$} ;
\draw (3,1.667) -- (3,5) node[above] {$e_2 - e_3$} ;
\draw (3,1.2) node {$e_1 - e_4$} ;
\draw (0,-0.2) node {$e_1 - e_2$} ;
\draw (6,-0.2) node {$e_3 - e_4$} ;
\draw (3,0.7) node {$\mathfrak{c}_7$} ;
\draw (4,1.5) node {$\mathfrak{c}_1$} ;
\draw (2,1.5) node {$\mathfrak{c}_2$} ;
\draw (3.5,2.2) node {$\mathfrak{c}_3$} ;
\draw (2.5,2.2) node {$\mathfrak{c}_4$} ;
\draw (3.5,3.4) node {$\mathfrak{c}_5$} ;
\draw (2.5,3.4) node {$\mathfrak{c}_6$} ;
\end{tikzpicture}
\end{center}
\begin{center}
\textup{Figure 2 \,:\, The chamber of $A_{3}^{+}$.}
\end{center}

\vspace{0.5cm}

\textup{For example, we set $m_{1,2} = 1$, $m_{1,3} = 1$, $m_{1,4} = 2$, $m_{2,3} = 1$, $m_{2,4} = 2$, and $m_{3,4} = 2$. For $a = \sum_{i=1}^{3} a_i(e_i - e_4) \in \overline{\mathfrak{c}_{{\rm nice}}}$,} 
\begin{eqnarray*}
P_{A_{3}^{+}, m}(a) = \left\{ (y_{i,j,k}) \in \mathbb{R}^{9} \,\middle|\, \begin{array}{ll} y_{i,j,k} \geq 0 \\ y_{1,2,1} + y_{1,3,1} + y_{1,4,1} + y_{1,4,2} = a_1 \\ - y_{1,2,1} + y_{2,3,1} + y_{2,4,1} + y_{2,4,2} = a_2 \\ - y_{1,3,1} - y_{2,3,1} + y_{3,4,1} + y_{3,4,2} = a_3 \end{array} \right\}. 
\end{eqnarray*}
\textup{From Corollary 2.10, we have}
\begin{align*}
v_{A_{3}^{+}, m, \mathfrak{c}_{\rm nice}}(a) &= {\rm Ires}_{x = 0}\left( \frac{e^{a_{1}x_1 + a_{2}x_2 + a_{3}x_3}}{x_{1}^{2}x_{2}^{2}x_{3}^2(x_1 - x_2)(x_1 - x_3)(x_2 - x_3)} \right) \\
&= \frac{1}{360} a_{1}^3(a_{1}^3 + 6a_{1}^{2}a_2 + 3a_{1}^{2}a_3 + 15a_{1}a_{2}^2 + 15a_{1}a_{2}a_3 + 10a_{2}^3 + 30a_{2}^{2}a_3).
\end{align*}
\end{exa}

\section{Main theorems}

In this section, we prove the main theorems of this paper. Let $\mathfrak{c}_{{\rm nice}}$ be the nice chamber of $A_{r}^{+}$ and let $a = \sum_{i = 1}^{r} a_{i}(e_i - e_{r+1}) \in \overline{\mathfrak{c}_{{\rm nice}}}$.

\begin{thm} 
For $a \in \overline{\mathfrak{c}_{{\rm nice}}}$, let $P_{A_{r}^{+}, m}(a)$ be the flow polytope as in Definition $2.1$ and let $v_{A_{r}^{+}, m, \mathfrak{c}_{\rm nice}}(a)$ be the volume of $P_{A_{r}^{+}, m}(a)$. Then $v = v_{A_{r}^{+}, m, \mathfrak{c}_{\rm nice}}(a)$ satisfies the system of differential equations as follows:
\begin{equation} 
\begin{cases}
\partial_{r}^{m_{r,r+1}} v = 0 \\
(\partial_{r-1} - \partial_{r})^{m_{r-1,r}} \partial_{r-1}^{m_{r-1,r+1}} v = 0 \\
\,\,\,\,\,\,\,\,\,\, \vdots  \\
(\partial_1 - \partial_2)^{m_{1,2}} (\partial_1 - \partial_3)^{m_{1,3}} \cdots (\partial_1 - \partial_r)^{m_{1,r}} \partial_{1}^{m_{1,r+1}} v = 0,
\end{cases}
\end{equation} 
where $\partial_i = \frac{\partial}{\partial a_i}$ for $i = 1, \dots, r$.
\end{thm}
\textit{Proof}. Let $F = \frac{e^{a_{1}x_1 + \cdots + a_{r}x_r}}{\prod_{i=1}^{r} x_{i}^{m_{i,r+1}} \prod_{1 \leq i < j \leq r} (x_i - x_j)^{m_{i,j}}}$. It is easy to see that
\begin{align*}
P(\partial_1, \dots, \partial_r)({\rm Ires}_{x=0} F) 
&= {\rm Ires}_{x=0}(P(\partial_1, \dots, \partial_r) F) \\
&= {\rm Ires}_{x=0}(P(x_1, \dots, x_r) F) 
\end{align*}
where $P$ is a polynomial.
Therefore, from Corollary 2.10, we obtain  
\begin{align*}
\partial_{r}^{m_{r,r+1}} v &= \partial_{r}^{m_{r,r+1}}{\rm Ires}_{x=0} F 
= {\rm Ires}_{x=0}\partial_{r}^{m_{r,r+1}} F \\
&= {\rm Ires}_{x=0}\left( \frac{e^{a_{1}x_1 + \cdots + a_{r}x_r}}{\prod_{i=1}^{r-1} x_{i}^{m_{i,r+1}} \prod_{1 \leq i < j \leq r} (x_i - x_j)^{m_{i,j}}} \right) \\
&= 0,
\end{align*}
and 
\begin{align*}
&(\partial_{r-1} - \partial_{r})^{m_{r-1,r}} \partial_{r-1}^{m_{r-1,r+1}} v \\
&= {\rm Ires}_{x=0} (\partial_{r-1} - \partial_{r})^{m_{r-1,r}} \partial_{r-1}^{m_{r-1,r+1}} F \\
&= {\rm Ires}_{x=0} (\partial_{r-1} - \partial_{r})^{m_{r-1,r}} \left( \frac{e^{a_{1}x_1 + \cdots + a_{r}x_r}}{x_{r}^{m_{r,r+1}} \prod_{i=1}^{r-2} x_{i}^{m_{i,r+1}} \prod_{1 \leq i < j \leq r} (x_i - x_j)^{m_{i,j}}} \right) \\
&= {\rm Ires}_{x=0} \left( \frac{e^{a_{1}x_1 + \cdots + a_{r}x_{r}}}{x_{r}^{m_{r,r+1}} \prod_{i=1}^{r-2} x_{i}^{m_{i,r+1}} \prod_{1 \leq i < j \leq r, (i,j) \neq (r-1,r)} (x_i - x_j)^{m_{i,j}}} \right) \\
&= {\rm Res}_{x_1 = 0} \cdots \left( {\rm Res}_{x_{r-1} = 0} \left( \frac{e^{a_{1}x_1 + \cdots + a_{r-1}x_{r-1}}}{\prod_{i=1}^{r-2} x_{i}^{m_{i, r+1}} \prod_{1 \leq i < j \leq r-1} (x_i - x_j)^{m_{i,j}}} \right. \right. \\
&\,\, \left. \left. \times {\rm Res}_{x_r = 0} \left( \frac{e^{a_{r}x_r}}{x_{r}^{m_{r,r+1}} \prod_{i = 1}^{r-2} (x_i - x_r)^{m_{i,r}}} \right) \right) \right)\\
&= 0,
\end{align*}
where we used
\[ {\rm Res}_{x_k = 0}\left( \frac{e^{a_{1}x_1 + \cdots + a_{k}x_k}}{\prod_{i=1}^{k-1} x_{i}^{m_{i,r+1}} \prod_{1 \leq i < j \leq k} (x_i - x_j)^{m_{i,j}}} \right) = 0 \]
for $k = 1, \dots, r$. 
Similarly, we can check the left expressions. \,\, $\square$

\begin{rema} 
{\rm In general, it is known that the volume function $v_{A_{r}^{+}, m}(a)$ of $P_{A_{r}^{+}, m}(a)$, as a distribution on $V$, satisfies the differential equation 
\[ Lv_{A_{r}^{+}, m}(a) = \delta(a), \] 
where $L = \prod_{i < j} (\partial_i - \partial_j)^{m_{i,j}}$ and $\delta(a)$ is the Dirac delta function on $V$ (\cite{6,9}). Note that $\partial_{r+1}$ in the definition of $L$ is supposed to be zero. The above theorem, together with Proposition 3.3 and Theorem 3.4 as below, characterizes the function $v_{A_{r}^{+}, m, \mathfrak{c}_{\rm nice}}(a)$ on $\overline{\mathfrak{c}_{\rm nice}}$ more explicitly.}
\end{rema}

Let $M_l = \sum_{i = l+1}^{r+1} m_{l,i} \,\, (l = 1, \dots, r)$. Then we have the following proposition.

\begin{propo} 
The coefficient of $a_{1}^{M_1 - 1}a_{2}^{M_2 - 1} \cdots a_{r-1}^{M_{r-1} - 1}a_{r}^{M_r - 1}$ in the volume function $v_{A_{r}^{+}, m}(a)$ is given by
\[ \frac{1}{(M_1 - 1) ! (M_2 - 1) ! \cdots (M_{r-1} - 1) ! (M_r - 1) !}. \]
\end{propo}
\textit{Proof}. 
From Corollary 2.10, we have 
\begin{equation*}
v_{A_{r}^{+}, m, \mathfrak{c}_{\rm nice}}(a) = \sum_{| i | = l - r} \frac{a_{1}^{i_1}}{i_{1} !}\frac{a_{2}^{i_2}}{i_{2} !} \cdots \frac{a_{r}^{i_r}}{i_{r} !} {\rm Ires}_{x = 0} \left( \frac{x_{1}^{i_1}x_{2}^{i_2} \cdots x_{r}^{i_r}}{\prod_{i = 1}^{r} x_{i}^{m_{i,r+1}} \prod_{1 \leq i < j \leq r} (x_i - x_j)^{m_{i,j}}} \right),
\end{equation*}
where $| i | = i_1 + \cdots + i_r$. 
When $i_l = M_l - 1$ for $l = 1, \dots, r$, 
\begin{align*}
&{\rm Ires}_{x = 0} \left( \frac{x_{1}^{M_1 - 1}x_{2}^{M_2 - 1} \cdots x_{r}^{M_r - 1}}{\prod_{i=1}^{r} x_{i}^{m_{i,r+1}} \prod_{1 \leq i < j \leq r} (x_i - x_j)^{m_{i,j}}} \right) \\
&= {\rm Res}_{x_1 = 0} \cdots {\rm Res}_{x_{r-1} = 0}{\rm Res}_{x_r = 0} \left( \frac{x_{1}^{(\sum_{i=2}^{r} m_{1,i}) - 1}x_{2}^{(\sum_{i=3}^{r} m_{2,i}) - 1} \cdots x_{r-1}^{m_{r-1,r} - 1}}{x_{r} \prod_{1 \leq i < j \leq r} (x_i - x_j)^{m_{i,j}}} \right) \\
&= {\rm Res}_{x_1 = 0} \cdots {\rm Res}_{x_{r-1} = 0} \left( \frac{x_{1}^{(\sum_{i=2}^{r-1} m_{1,i}) - 1}x_{2}^{(\sum_{i=3}^{r-1} m_{2,i}) - 1} \cdots x_{r-2}^{m_{r-2,r-1}-1}}{x_{r-1} \prod_{1 \leq i < j \leq r-1} (x_i - x_j)^{m_{i,j}}} \right) \\
&= {\rm Res}_{x_1 = 0} \frac{1}{x_1} = 1.
\end{align*}
Thus we obtain the proposition. \,\, $\square$

\begin{thm} 
Let $\phi_r = \phi(a_1, \dots, a_r)$ be a homogeneous polynomial of $a_1, \dots, a_r$ with degree $d$ and let $M = \sum_{1 \leq i < j \leq r+1} m_{i,j}$. Suppose $\phi_r$ satisfies the system of differential equations as follows:
\begin{equation} 
\begin{cases}
\partial_{r}^{m_{r,r+1}}\phi_r = 0 \\
(\partial_{r-1} - \partial_{r})^{m_{r-1,r}}\partial_{r-1}^{m_{r-1,r+1}}\phi_r = 0 \\
\,\,\,\,\,\,\,\,\,\, \vdots  \\
(\partial_1 - \partial_2)^{m_{1,2}}(\partial_1 - \partial_3)^{m_{1,3}} \cdots (\partial_1 - \partial_r)^{m_{1,r}}\partial_{1}^{m_{1,r+1}}\phi_r = 0. 
\end{cases}
\end{equation}
\begin{itemize}
\item[$(1)$] If $M - r < d$, then $\phi_r = 0$.
\item[$(2)$] If $0 \leq d \leq M - r$, then there is a non trivial homogeneous polynomial $\phi_r$ satisfying $(3.2)$.
\item[$(3)$] If $d = M - r$ in particular, $\phi_r$ is equal to a constant multiple of $v = v_{A_{r}^{+}, m, \mathfrak{c}_{\rm nice}}(a)$.
\end{itemize}
\end{thm}
\textit{Proof}. We argue by induction on $r$. In the case that $r = 1$, we write
\[ \phi_1 = \phi(a_1) = p a_{1}^d, \]
where $p$ is a constant. 
If $m_{1,2} - 1 < d$ and $\phi_1$ satisfies the differential equation $\partial_{1}^{m_{1,2}} \phi_1 = 0$, then $p = 0$ and hence $\phi_1 = 0$. If $0 \leq d \leq m_{1,2} - 1$, then for any $p \neq 0$, $\partial_{1}^{m_{1,2}} \phi_1 = 0$.
Also, if $d = m_{1,2} - 1$, in particular, then $\phi_1 = pa_{1}^{m_{1,2} - 1}$, while $v = \frac{1}{(m_{1,2} - 1) !}a_{1}^{m_{1,2} - 1}$ as in Example 2.11. Hence $\phi_1$ is equal to a constant multiple of $v$.

We assume that the statement of this theorem holds for $r-1$. We write $\phi_r$ as 
\[ \phi_r = \phi(a_1, \dots, a_r) = g_d(a_2, \dots, a_r) + a_{1}g_{d-1}(a_2, \dots, a_r) + \cdots + a_{1}^{d}g_0(a_2, \dots, a_r), \]
where $g_k$ is a homogeneous polynomial of $a_2, \dots, a_r$ with degree $k$ for $k = 0, 1, \dots, d$. Then for $k = 0, 1, \dots, d$, $g_k$ satisfies the differential equations as follows:
\begin{equation} 
\begin{cases}
\partial_{r}^{m_{r,r+1}} g_k = 0 \\
(\partial_{r-1} - \partial_{r})^{m_{r-1,r}}\partial_{r-1}^{m_{r-1,r+1}} g_k = 0 \\
\,\,\,\,\,\,\,\,\,\, \vdots \\
(\partial_2 - \partial_3)^{m_{2,3}}(\partial_2 - \partial_4)^{m_{2,4}} \cdots (\partial_2 - \partial_{r})^{m_{2,r}}\partial_{2}^{m_{2,r+1}}g_k = 0. 
\end{cases}
\end{equation}
We set $h = (\sum_{2 \leq i < j \leq r+1} m_{i,j}) - (r - 1)$. From the inductive assumption, if $0 \leq k \leq h$, then $g_k$ is a homogeneous polynomial. On the other hand, if $h+1 \leq k \leq d$, then $g_k = 0$, namely,
\begin{equation} 
g_d(a_2, \dots, a_r) = g_{d-1}(a_2, \dots, a_r) = \cdots = g_{h + 1}(a_2, \dots, a_r) = 0.
\end{equation}

\hspace{-0.5cm} (1) We consider the case of $M - r < d$. Let $M_{1} = \sum_{i=2}^{r+1} m_{1,i}$. 
Now we compare the coefficients of $a_{1}^{d - h - M_1 + n}$ in $(\partial_1 - \partial_2)^{m_{1,2}}(\partial_1 - \partial_3)^{m_{1,3}} \cdots (\partial_1 - \partial_{r})^{m_{1,r}}\partial_{1}^{m_{1,r+1}}\phi_r$ for $n = 0, \dots, h$. For $q = 1, \dots, M_1 - m_{1,r+1}$, we define 
\begin{align} 
D_q &= \displaystyle \sum_{2 \leq i_1 \leq r} \left( \begin{smallmatrix} m_{1, i_1} \\ q \end{smallmatrix} \right) \partial_{i_1}^{q} + \cdots + \sum_{\substack{p_1 + \cdots + p_k = q \\ 2 \leq i_1 < \cdots < i_k \leq r}} \left( \prod_{1 \leq l \leq k} \left( \begin{smallmatrix} m_{1,i_l} \\ p_l \end{smallmatrix} \right) \right) \partial_{i_1}^{p_1}\partial_{i_2}^{p_2} \cdots \partial_{i_k}^{p_k} \notag \\
&+ \cdots + \sum_{2 \leq i_1 < \cdots < i_q \leq r} \left( \prod_{1 \leq l \leq q} \left( \begin{smallmatrix} m_{1,i_l} \\ 1 \end{smallmatrix} \right) \right) \partial_{i_1}\partial_{i_2} \cdots \partial_{i_q}.
\end{align}
Then we have the following equation:
\begin{align} 
&\frac{(d - h + n) !}{(d - h - M_1 + n) !}g_{h - n}(a_2, \dots, a_r) - \frac{(d - h + n - 1) !}{(d - h - M_1 + n) !}D_1 g_{h - n + 1}(a_2, \dots, a_r) \notag \\ 
&+ \cdots \pm \frac{(d - h + n - j) !}{(d - h - M_1 + n) !}D_j g_{h - n + j}(a_2, \dots, a_r) \notag \\
&\pm \cdots \pm \frac{(d - h + n - (M_1 - m_{1,r+1})) !}{(d - h - M_1 + n) !} D_{M_1 - m_{1,r+1}} g_{h - n + M_{1,r}}(a_2, \dots, a_r) \notag \\
&= 0.
\end{align}
When $n = 0$, from (3.4) and (3.6), we have
\begin{equation*}
g_{h}(a_2, \dots, a_r) = 0.
\end{equation*}
When $n = 1$, we have
\[ \frac{(d - h + 1) !}{(d - h - M_1 + 1) !}g_{h - 1}(a_2, \dots, a_r) - \frac{(d - h) !}{(d - h - M_1 + 1) !}D_1 g_{h}(a_2, \dots, a_r) = 0. \] 
Thus we have
\[ g_{h - 1}(a_2, \dots, a_r) = 0. \]
Similarly, we have
\[ g_{h - 2}(a_2, \dots, a_r) = g_{h - 3}(a_2, \dots, a_r) = \cdots = g_0(a_2, \dots, a_r) = 0 \]
and hence $\phi_r = 0$. 

\hspace{-0.5cm} (2) We consider the case of $0 \leq d \leq M - r$. By the inductive assumption, there is a non trivial homogeneous polynomial $g_{h - n_1 + i}$ satisfying (3.3) for $i = 1, \dots, n_1$, where $n_1 = M - r - d + 1$. We can take 
\[ g_{h - n_1 + i}(a_2, \dots, a_r) \neq 0. \]
When $n = n_1$, from (3.4) and (3.6),
\begin{align*}
&g_{h - n_1}(a_2, \dots, a_r) \\
&= \frac{(d - h + n_1 - 1) !}{(d - h + n_1) !}D_1 g_{h - n_1 + 1}(a_2, \dots, a_r) - \frac{(d - h + n_1 - 2) !}{(d - h + n_1) !}D_2 g_{h - n_1 + 2}(a_2, \dots, a_r) \\ 
&+ \cdots \pm \frac{(d - h) !}{(d - h + n_1) !}D_{n_1} g_{h}(a_2, \dots, a_r).
\end{align*}
When $n = n_1 + 1$,
\begin{align*} 
&g_{h - (n_1 + 1)}(a_2, \dots, a_r) \\
&= \frac{(d - h + n_1) !}{(d - h + n_1 + 1) !}D_1 g_{h - n_1}(a_2, \dots, a_r) - \frac{(d - h + n_1 - 1) !}{(d - h + n_1 + 1) !}D_2 g_{h - n_1 + 1}(a_2, \dots, a_r) \\ 
&+ \cdots \pm \frac{(d - h) !}{(d - h + n_1 + 1) !}D_{n_1 + 1} g_{h}(a_2, \dots, a_r).
\end{align*}
Similarly, 
for $n = n_1 + 2, \dots, h$, we can express $g_{h - j}(a_2, \dots, a_r) \,\, (j = n_1, n_1 + 1, \dots, h)$ in terms of $g_{h - j + i}(a_2, \dots, a_r) \,\, (i = 1, \dots, j)$ and their partial derivatives. Namely, we can express $\phi_r$ in terms of $g_{h - n_1 + i}(a_2, \dots, a_r)$ and their partial derivatives. It follows that $\phi_r \neq 0$ when $0 \leq d \leq M - r$. 

\hspace{-0.5cm} (3) If $d = M - r$ in particular, then $n_1 = 1$, and $g_{h - j} \,\, (j = 1, \dots, h)$ becomes linear combination of $g_h$ and their partial derivatives. Therefore $\phi_r$ is  uniquely determined by $g_{h}$. Moreover, from the inductive assumption, $g_h = C \cdot v_{A_{r-1}^{+}, m', \mathfrak{c}'_{\rm nice}}$, where $C$ is a constant, $m' = (m_{i,j})_{2 \leq i < j \leq r+1}$, and $\mathfrak{c}'_{\rm nice}$ is a nice chamber of $A_{r-1}^{+}$. Hence the solution of (3.2) is unique up to a constant multiple. On the other hand, by Theorem 3.1, $v_{A_{r}^{+}, m, \mathfrak{c}_{\rm nice}}$ satisfies the system of differential equations (3.2). Hence $\phi_r$ is equal to a constant multiple of $v_{A_{r}^{+}, m, \mathfrak{c}_{\rm nice}}$. \,\, $\square$

\begin{rema} 
\textup{Let $M_1 = \sum_{i=2}^{r+1} m_{1,i}$ and let $D_q \,\, (q = 1, \dots, h)$ be as in $(3.5)$. When $d = M - r$, from the proof of Theorem $3.4 \, (3)$, $g_{h-j} \,\, (j = 1, \dots, h)$ is uniquely determined as follows:
\begin{equation*}
\begin{cases}
g_{h-1} = \frac{(M_1 - 1) !}{M_1 !} D_1 g_{h} \\ 
g_{h-2} = \frac{(M_1 - 1) !}{(M_1 + 1) !} (D_{1}^2 - D_2) g_{h} \\ 
g_{h-3} = \frac{(M_1 - 1) !}{(M_1 + 2) !} (D_{1}^3 - 2D_{1}D_2 + D_3) g_h \\
\,\,\,\,\,\,\,\,\,\, \vdots \\ 
g_{0} = \frac{(M_1 - 1) !}{(M - r) !} (D_{1}^{h} - (h - 1)D_{1}^{h-2}D_2 + \cdots \pm D_{h}) g_{h}.
\end{cases}
\end{equation*}}
\end{rema}


Let $m' = (m_{i,j})_{2 \leq i < j \leq r+1}$, $\mathfrak{c}'_{{\rm nice}}$ a nice chamber of $A_{r-1}^{+}$ and $a' = \sum_{i = 2}^{r} a_{i}(e_i - e_{r+1}) \in \overline{\mathfrak{c}'_{{\rm nice}}}$. From Proposition 3.3 and Remark 3.5, we obtain the following theorem.

\begin{thm} 
Let $h = (\sum_{2 \leq i < j \leq r+1} m_{i,j}) - (r - 1)$ and let $D_q \,\, (q = 1, \dots, h)$ be as in $(3.5)$. Then $v_{A_{r}^{+}, m, \mathfrak{c}_{\rm nice}}(a)$ is written by linear combination of $v_{A_{r-1}^{+}, m', \mathfrak{c}'_{\rm nice}}(a')$ and its partial derivatives as follows:
\begin{align}
v_{A_{r}^{+}, m, \mathfrak{c}_{\rm nice}}(a)  
&= \left\{ \frac{1}{(M_1 - 1) !}a_{1}^{M_1 - 1} + \frac{1}{M_1 !} a_{1}^{M_1} D_1 + \frac{1}{(M_1 + 1) !} a_{1}^{M_1 + 1} (D_{1}^2 - D_2) \right. \notag \\
&\left. + \frac{1}{(M_1 + 2) !} a_{1}^{M_1 + 2} (D_{1}^3 - 2D_{1}D_2 + D_3) + \cdots \right. \notag \\ 
&\left. + \frac{1}{(M - r) !} a_{1}^{M - r} (D_{1}^h - (h - 1)D_{1}^{h-1}D_2 + \cdots \pm D_h) \right\} v_{A_{r-1}^{+}, m', \mathfrak{c}'_{\rm nice}}(a'). \,\,\,\,\, 
\end{align} 
\end{thm}

\begin{exa} 
\textup{Let $r = 3$, let $a = \sum_{i = 1}^{3} a_{i}(e_i - e_{4}) \in \overline{\mathfrak{c}_{{\rm nice}}}$ and let $a' = \sum_{i = 2}^{3} a_{i}(e_i - e_{4}) \in \overline{\mathfrak{c}'_{{\rm nice}}}$. We set $m_{1,2} = 1, m_{1,3} = 1, m_{1,4} = 2, m_{2,3} = 1, m_{2,4} = 2$ and $m_{3,4} = 2$ as in Example 2.13. Then we have}
\[ v_{A_{3}^{+}, m, \mathfrak{c}_{\rm nice}}(a) = \frac{1}{360} a_{1}^3(a_{1}^3 + 6a_{1}^{2}a_2 + 3a_{1}^{2}a_3 + 15a_{1}a_{2}^2 + 15a_{1}a_{2}a_3 + 10a_{2}^3 + 30a_{2}^{2}a_3). \]
\textup{We can check that $v = v_{A_{3}^{+}, m, \mathfrak{c}_{\rm nice}}(a)$ satisfies the system of differential equations as follows:}
\begin{equation*}
\begin{cases}
\partial_{3}^{2}v = 0 \\
(\partial_2 - \partial_3)\partial_{2}^{2}v = 0 \\
(\partial_1 - \partial_2)(\partial_1 - \partial_3)\partial_{1}^{2}v = 0.
\end{cases}
\end{equation*}
\textup{Also, from Proposition 3.2, the coefficient of the term $a_{1}^{3}a_{2}^{2}a_3$ is $\frac{1}{3 ! 2 ! 1!} = \frac{1}{12}$. When $r = 2$,}
\[ v_{A_{2}^{+}, m', \mathfrak{c}'_{\rm nice}}(a') = \frac{1}{6}a_{2}^2(a_2 + 3a_3). \]
\textup{Therefore, we have}
\begin{align*}
&\left\{ \frac{1}{6}a_{1}^3 + \frac{1}{24}a_{1}^{4}D_1 + \frac{1}{120}a_{1}^{5}(D_{1}^2 - D_2) + \frac{1}{720}a_{1}^{6}(D_{1}^3 - 2D_{1}D_2 + D_3) \right\}v_{A_{2}^{+}, m', \mathfrak{c}'_{\rm nice}}(a') \\
&= \frac{1}{36}a_{1}^{3}a_{2}^{3} + \frac{1}{12}a_{1}^{3}a_{2}^{2}a_3 + \frac{1}{24}a_{1}^{4}a_{2}^2 + \frac{1}{24}a_{1}^{4}a_{2}a_3 + \frac{1}{60}a_{1}^{5}a_2 + \frac{1}{120}a_{1}^{5}a_3 + \frac{1}{360}a_{1}^6 \\
&= v_{A_{3}^{+}, m, \mathfrak{c}_{\rm nice}}(a).
\end{align*}
\textup{Hence when $r = 3$, we can check the equation (3.7) in Theorem 3.6.}
\end{exa}

\vspace{0.5cm}

\hspace{-0.75cm} \textbf{Acknowledgements} \\

\hspace{-0.75cm} The third author was supported by JSPS KAKENHI Grant Number JP16K05137.



\hspace{-0.75cm} {\footnotesize Takayuki NEGISHI \\
Department of Mathematics, Graduate School of Science and Engineering, Chuo University, \\
Kasuga, Bunkyo-Ku, Tokyo, 112--8551 Japan. \\
\vspace{0.5mm} \\
Yuki SUGIYAMA \\
Department of Mathematics, Graduate School of Science and Engineering, Chuo University, \\
Kasuga, Bunkyo-Ku, Tokyo, 112--8551 Japan. \\
\textit{e-mail}: y-sugi@gug.math.chuo-u.ac.jp \\
\vspace{0.5mm} \\
Tatsuru TAKAKURA \\
Department of Mathematics, Chuo University, \\
Kasuga, Bunkyo-Ku, Tokyo, 112--8551 Japan. \\
\textit{e-mail}: takakura@math.chuo-u.ac.jp}

\end{document}